\documentclass[12pt]{article}
\usepackage{amssymb,amsmath,amsfonts,amssymb}
\usepackage{graphicx}
\input{epsf.sty}
\newtheorem{theorem}{\NI{\bf Theorem}}[section]
\newtheorem{lemma}{\NI\bf Lemma}[section]

\newtheorem{defn}{\NI\bf Definition}[section]

\newtheorem{remark}{\NI\bf Remark}[section]
\newcommand{\bt}{\begin{theorem}}
\newcommand{\et}{\end{theorem}}
\newcommand{\NI}{\noindent}
\newcommand{\newsection}[1]{\setcounter{equation}{0}
\setcounter{theorem}{0}
\section{#1}}

\newcommand{\bea}{\begin{eqnarray}}
\newcommand{\eea}{\end{eqnarray}}
\newcommand{\ben}{\begin{eqnarray*}}
\newcommand{\een}{\end{eqnarray*}}

\newcommand {\olf}{\overline {F}}

\newcommand {\If}{F^{-1}}
\newcommand {\Ig}{G^{-1}}
\newcommand {\IR}{I\!\!R}

\newcommand{\lmd}{\lambda}

\newcommand{\Raro}{\Rightarrow}
\newcommand{\Laro}{\Leftrightarrow}

\newcommand{\vsp}{\vskip 1em}

\def \qed {\hfill \vrule height 6pt width 6pt depth 0pt}
\newcommand{\be}{\begin{equation}}
\newcommand{\ee}{\end{equation}}

\title {Dependence comparisons of order  statistics in the proportional hazards model }
\author{Subhash Kochar\\
 Fariborz Maseeh Department of Mathematics and Statistics\\
Portland State University, Portland, OR, USA\\
{Email:kochar@pdx.edu}}

\begin {document}
\maketitle
\vsp

\begin{abstract}
Let $X_1, \ldots , X_n$ be mutually independent exponential random variables with distinct hazard rates $\lambda_1, \ldots , \lambda_n > 0$ and let $Y_1, \ldots, Y_n$  be a random sample from the exponential distribution with hazard rate $\bar \lmd = \sum_{i=1}^n \lmd_i/n$.  Also let $X_{1:n} < \cdots < X_{n:n}$ and $Y_{1:n} < \cdots < Y_{n:n}$ be their associated order statistics. It is shown that for $1\le i <j \le n$, the generalized spacing $X_{j:\, n} - X_{i:\, n}$ is more dispersed than $Y_{j:\,n} - Y_{i:\, n}$ according to dispersive ordering. This result is used to solve  a long standing open problem that for $2\le i \le n$ the dependence of 
$ X_{i:\, n}$ on  $X_{1:\, n}$  is less than that of $Y_{i: \, n}$ on $Y_{1\, :n}$, in the sense of the  more stochastically increasing. This dependence  result is also extended to  the PHR model. This extends  the  earlier work of  {\em Genest, Kochar and Xu}[   J.\ Multivariate Anal.\  {\bf 100} (2009) \ 1587-1592] who proved this result for  $i =n$.

\end{abstract}
\newsection{Introduction}
\vsp
Several  notions of monotone dependence have been introduced and studied in the literature. Researchers have also  developed the corresponding dependence (partial) orderings which compare the degree of (monotone) dependence within the components of different random  vectors of the same length. For details, see the pioneering paper of {\em Lehmann} (1966) and Chapter 5 of {\em Barlow and  Porschan }(1981)  for  different notions of positive dependence, and that of {\em Kimeldorf and Sampson} (1989) for a unified treatment of families, orderings and measures of monotone dependence.  For more detailed   discussion of these concepts see   Chapter 2 of {\em Joe} (1997) and Chapter 5 of {\em  Nelsen} (1999).

Let $\{X_1, \ldots, X_n\}$ be a set of  continuous  random variables. Many authors have investigated the nature of the dependence that may exist between two order statistics $X_{i:n}$ and $X_{j:n}$ for $1 \le i < j \le n$ under different distributional scenarios. It seems natural to expect some degree of positive dependence between them. When $X_1,\ldots,X_n$ are independent (but not necessarily identically distributed) the order statistics $X_{1:n},\ldots,X_{n:n}$ are associated since the order statistics are increasing functions of the sample observations and independent random variables are associated.  This yields many useful product inequalities for order statistics of independent random variables, and in particular, $ Cov(X_{i:n},X_{j:n}) \geq 0$ for all $i$ and $j$, which was initially proved by {\em Bickel} (1967) when the $X$'s are  independent and identically distributed (i.i.d.). {\em Boland et al.} (1996) studied this problem in detail and discussed  different types of dependence that hold between them. It is shown in that paper that in the case of i.i.d. observations, any pair of order statistics is  $TP_2$ dependent  (also called monotone likelihood ratio dependent) which is the strongest type of dependence in the hierarchy of various dependence criteria as described in {\em Barlow and Proschan} (1981). But this is not the case in  the non-iid case. However, as  is shown in {\em Boland et al.} (1996) that in, general, whereas  for $ i < j$, $X_{j:n}$ may not be stochastically increasing in $X_{i:n}$, \  $X_{j:n}$ is always right tail increasing (RTI) in $X_{i:n}$. The reader is  referred to Chapter 5 of {\em Barlow and Proschan} (1981) for basic definitions and relations among various types of dependence.

It is also of interest to compare the strength of dependence  that may exist between two pairs of order statistics. When the parent distribution from which the random sample is drawn  has an increasing hazard rate and a decreasing reverse hazard rate, {\em Tukey} (1958) showed that
\begin{equation}
\label{ine1}
{\rm Cov}(X_{i':n}, X_{j':n})\le {\rm Cov}(X_{i:n}, X_{j:n})
\end{equation}
for either $i=i'$ and $j\le j'$; or $j=j'$ and $i'\le i$. {\em Kim and David} (1990) proved that if both the hazard  and the reverse
hazard rates of the $X_i$'s are increasing, then  inequality \eqref{ine1}  remains valid when  $i=i'$
and $j\le j'$; However, the inequality \eqref{ine1} is reversed when $j=j'$ and $i'\le i$.

Let $X_1, \ldots , X_n$ are mutually independent exponentials with distinct hazard rates $\lambda_1, \ldots , \lambda_n > 0$ and let  $Y_1, \ldots, Y_n$  form a random sample from the exponential distribution with hazard rate $\overline {\lmd} = ( \sum_{i=1}^n X_i)/n$.
 {\em Sathe} (1988)  proved  that for any $j \in \{2,\ldots,n\}$
\be
\label{Sathe}
\mbox{corr} (X_{1:n}, X_{j:n} )\le {corr} (Y_{1:n}, Y_{j:n} ).
\ee
Although this observation is interesting, it merely compares the relative degree of \emph{linear} association within the two pairs. It is now widely recognized, however, that margin-free measures of association are more appropriate than Pearson's correlation, because they are based on the unique underlying copula which governs the dependence between the components of a continuous random pair.

Section 2 is on preliminaries where various definitions and notations are given.  The main results of this paper are given in Chapter 3.
\vsp
\newsection{Preliminaries}

For $i = 1, 2$, let $(X_i, Y_i)$ be a pair of continuous random
variables with joint cumulative function $H_i$ and margins $F_i$,
$G_i$. Let
$$
C_i (u,v) = H_i \{ F_i^{-1}(u), G_i^{-1}(v) \}, \quad u,v \in (0,1)
$$
be the unique copula associated with $H_i$. In other words, $C_i$ is
the distribution of the pair $(U_i,V_i) \equiv (F_i(X_i),G_i(Y_i))$
whose margins are uniform on the interval $(0,1)$. See, e.g.,
Chapter~1 of {\em Nelsen} (1999) for details.

The most well understood partial order to compare the strength of dependence within two random vectors is that of 
{\em positive quadrant dependence} (PQD) as defined  below.
\begin{defn}
A copula $C_1$ is said to be less dependent than copula $C_2$ in the \emph{positive quadrant dependence ordering} (PQD), denoted $(X_1,
Y_1) \prec_{\rm PQD} (X_2, Y_2)$, if and only if 
$$
C_1 (u,v) \le C_2 (u,v), \quad u, v \in (0,1).
$$
\end{defn}
This condition implies that $ \kappa (S_1, T_1) \le \kappa (S_2, T_2) $ for all concordance measures meeting the axioms of Scarsini
(1984) like Kendall's $\tau$ and Spearman's $\rho$.  See  Tchen (1980) for details.

{\em Lehmann} (1966) in his seminal work introduced the notion of  {\em  monotone regression dependence (MRD)} which is  also  known as   {\em  stochastic increasingness} ({\bf SI}) in the literature.

\begin{defn}
Let $(X,Y)$ be a  bivariate random vector with joint distribution function $H$.\  $Y$ is said to be stochastically  increasing (SI) in $X$ if for all $(x_1,x_2) \in \IR^2$,
\be \label{SI}
 x_1 < x_2 \Raro P( Y \le  y|\; X = x_2)  \le  P( Y \le y|\; X = x_1), \mbox{ for  all $y\in \IR$,}
\ee
\end{defn}
If we denote by $H_x$ the distribution function   of the conditional distribution of $Y$ given $X=x$, then \eqref{SI}  can be rewritten  as
\be
\label{SI-2}
x_1 < x_2\Raro H_{x_2}\circ H_{x_1}^{-1}(u) \le  u, \quad \mbox{for $0 \le u \le 1$}.
\ee

 Note that in case $X$ and $Y$ are  independent, $H_{x_2}\circ H_{x_1}^{-1}(u) = u$, for $0 \le u \le 1$ and for all $(x_1, x_2)$.
The  SI property is not symmetric in $X$ and $Y$. It is a very strong  notion of positive
dependence and many of the other notions of positive dependence follow from it. 

Denoting  by  $\xi_p = \If(p)$, the  $pth$ quantile of the marginal distribution of $X$, we see that \eqref{SI-2} will  hold if and only if  for all $0\le u \le 1$,
\be \label{copula-SI}
0 \le p < q \le 1 \Raro H_{\xi_q}\circ H _{{\xi_p}}^{-1}(u) \le  u.
\ee

 In his survey, {\em Joe} (1997) mentions a number of bivariate stochastic ordering relations. One such notion 
  is that of greater monotone regression (or more SI) dependence, originally considered by {\em Yanagimoto and Okamoto} (1969) 
  and later extended and further investigated by {\em Schriever} (1987),  {\em Cap\'era\`a and  Genest } (1990),  
    and {\em Av\'{e}rous, Genest and Kochar} (2005), among others. As discussed  in the books by 
    {\em Joe} (1997) and {\em Nelsen} (1999)  a reasonable    way  to compare  the relative degree of dependence 
     between two random vectors is through their copulas. We will discuss here the concept of  more SI, a partial order 
     to compare the strength of dependence that may exist between two bivariate random vectors in the sense of monotone 
     regression dependence (stochastic increasingness).

Suppose we have two pairs of continuous random variables $(X_1,Y_1)$   and $(X_2,Y_2)$   with joint cumulative distribution functions $H_i$ and marginals $F_i$ and $G_i$ for $i=1,2$.
\begin{defn}
\label{more-si-dfn} $Y_2$ is said to be more stochastically increasing in $X_2$ than $Y_1$ is in $X_1$, denoted by $(Y_1 | \,
X_1) \prec_{\rm SI} (Y_2 | \, X_2)$ or $H_1 \prec_{\rm SI} H_2$, if
\begin{equation}
\label{more-si} 0 < p \le q < 1 \Longrightarrow H_{2,\, \xi_{2q} }\circ{H^{-1}_{2,\, \xi_{2p}}}(u) \le H_{1,\, \xi_{1q}}\circ{H^{-1}_{1,\,\xi_{1p}}}(u),
\end{equation}
for all $u \in (0,1)$, where for $i=1,2$, $H_{i, s}$ denotes the conditional distribution of $Y_i$ given $X_i = s$, and $\xi_{ip} =
F_i^{-1} (p)$ stands for the $p$th quantile of the marginal distribution of  $X_i$.
\end{defn}

Obviously, $(\ref{more-si})$ implies  that $Y_2$ is stochastically increasing in $X_2$ if $X_1$ and $Y_1$ are independent. It also implies that if $Y_1$ is stochastically increasing in $X_1$, then so is $Y_2$ in $X_2$; and conversely, if $Y_2$ is stochastically decreasing in $X_2$, then so is $Y_1$ in $X_1$.The above definition of more SI is copula based and $$H_1 \prec_{\rm SI} H_2 \Laro H_1 \prec_{\rm PQD} H_2$$

 {\em Av\'{e}rous, Genest and Kochar} (2005) have shown in their paper that in the case of i.i.d. observations, the copula of any pair of order statistics is independent of the distribution of the parent observations as long as they are continuous. In a sense, their copula has a distribution-free property. But this is not the case  if the observations are not i.i.d.

The natural question is to see if  we can   extend \eqref{Sathe}  result to a copula based positive  dependence ordering.  {\em Genest, Kochar and Xu} (2009) proved that
 under the given conditions,
 
 \be
 \label{dep1}
(X_{n:n}|X_{1:n})\prec_{\rm SI}(Y_{n:n}|Y_{1:n}).
\ee
It has been an open problem to see whether this result can be generalized  from the the largest order statistics to other order statistics. We prove in this paper that  for $2\le i \le n$,
 $$
(X_{i:n}|X_{1:n})\prec_{\rm SI}(Y_{i:n}|Y_{1:n}).
$$

This implies in particular that
$$
\kappa (X_{i:n}, X_{1:n} ) \le \kappa (Y_{i:n}, Y{i:n} ),
$$
where $\kappa (S,T)$ represents any concordance measure between random variables $S$ and $T$ in the sense of Scarsini (1984), e.g., Spearman's rho or Kendall's tau.
A related work to this problem is {\em Dolati, Genest and Kochar} (2008).

\newsection {Main results}
\vsp

The proof of our main result   relies heavily on the notion of dispersive ordering between two random variables $X$ and $Y$, and properties thereof. For completeness, the definition of this
concept is recalled below.

\begin{defn}
\label{dispersive-definition} A random variable $X$ with distribution function $F$ is said to be less dispersed than another variable $Y$ with distribution $G$, written as $X \le_{disp} Y$
 or $F \le_{disp} G$, if and only if
$$
F^{-1} (v) - F^{-1} (u) \le G^{-1} (u) - G^{-1} (v)
$$
for all $0 < u \le v < 1$.
\end{defn}

Dispersive ordering  is closely related to star-order which is a partial order to compare the relative aging or skewness and is defined as below.
\begin{defn}
A random variable $X$ with distribution function $F$ is said to be star ordered   with respect to another random variable $Y$ with distribution $G$, written as $X \le_{*} Y$
 or $F \le_{*} G$, if and only if
 $$ \frac{\Ig \circ F(x)}{x} \,\, \mbox{ is increasing in $x$}.$$
 \end{defn}
For nonnegative random variables, star order is related to dispersive order by
$$X  \le_{*} Y \Leftrightarrow \log{ X} \le_{disp}\log{ Y}$$

The proof of the next lemma, which is a refined version of  a result by {\em Deshpande and Kochar} (1983), on relation between star-order and dispersive order, can be found in {\em Kochar and Xu} (2012).
\begin{lemma}
\label{disp-*}
Let  $X$ and $Y$ be two random variables. If $ X \le_{*} Y$ and $X \le_{st} Y$, then $X\le_{disp} Y$.
\end{lemma}



The following result of {\em Khaledi and Kochar } (2005) will be used to prove our main theorem.

\begin{lemma}
\label{disp-si}
({\em Khaledi and Kochar, 2005}) \, Let  $X_i$ and $Y_i$ be independent random variables with distribution functions $F_i$ and $G_i$, respectively for  $i=1,2$. Then
\be \label{SI-disp}
X_2 \le_{disp} X_1 \mbox{ and $ Y_1 \le_{disp} Y_2$} \Raro (X_2 +Y_2)| X_2 \prec_{\rm SI}(X_1  +Y_1)| X_1
\ee
\end{lemma}

The next theorem on dispersive ordering between general spacings  which is also of independent interest, will be  used  to  prove our main result.

\bt
\label{gen=space-disp}
Let $X_1, \ldots , X_n$ be mutually independent exponential random variables with distinct hazard rates $\lambda_1, \ldots , \lambda_n > 0$ and let $Y_1, \ldots, Y_n$  be a random sample from the exponential distribution with hazard rate $\bar \lmd$.  Then for $1 \le i < j \le n$,
\be
\label{Gspace-disp}
 (Y_{j:n} - Y_{i:n}) \le_{disp} (X_{j:n} - X_{i:n})
\ee
\et

\NI{\bf Proof.}\, {\em Yu} (2021) proved in Corollary 1 of his paper that for  $1 \le i < j \le n$,
 \be
\label{Gspace-*}
(Y_{j:n} - Y_{i:n}) \le_{*} (X_{j:n} - X_{i:n})
 \ee
and   {\em Kochar and Rojo} (1996) proved in  their Corollary 2.1 that
\be
\label{Gspace-st}
(Y_{j:n} - Y_{i:n}) \le_{st} (X_{j:n} - X_{i:n})
 \ee

Using  \eqref{Gspace-*} and \eqref{Gspace-st}, we get \eqref{Gspace-disp} using  Lemma \ref{disp-*}.

\qed

\vsp

Now we give the main result of this paper.
\bt
\label{dep2}
Let $X_1, \ldots , X_n$ are mutually independent exponentials with distinct hazard rates $\lambda_1, \ldots , \lambda_n > 0$ and let  $Y_1, \ldots, Y_n$  form a random sample from the exponential distribution with hazard rate $\overline {\lmd} = ( \sum_{i=1}^n \lmd_i)/n$. Then for $i \in \{2,\ldots,n\}$,
\be
\label{eqdep2}
(X_{i:n}|X_{1:n})\prec_{\rm SI}(Y_{i:n}|Y_{1:n}).
\ee
\et

\NI {\bf Proof.} \,\,  It follows from Theorem \ref{gen=space-disp} above that for $2\le i \le n$,
\be
\label{range-disp}
 (Y_{i:\,n} - Y_{1:\,n}) \le_{disp}(X_{i:\, n} -  X_{1:\, n}).
\ee

{\em Kochar and Korwar} (1996)  proved that $ (X_{i:n} - X_{1:n})  $ is independent of $X_{1:\,n}$ and $ X_{1:\,n} \stackrel{st}{=} Y_{1:\,n}$.   Similarly, $ (Y_{i: \, n} -  Y_{1: \, n})$ is independent of $Y_{1:\,n}$. Now we can express $X_{i:n}$ and  $Y_{i:n}$ as
  $$
  X_{i:n} = (X_{i: \, n} -  X_{1: \, n}) + X_{1:n} \mbox { and } Y_{i:\,n} =   (Y_{i:n} - Y_{1:n}) + Y_{1:n}
  $$
  
  Using \eqref{range-disp} and the above two two equations, it follows from Lemma \ref{disp-si} that
   \be
(X_{i:n}|X_{1:n})\prec_{\rm SI}(Y_{i:n}|Y_{1:n}), \, \mbox{ for $2\le i \le n$}. \ee

This proves the required result.
\qed
 \vsp

\begin{remark}.\, {\rm Since the copula of the order statistics of a random sample has the distribution-free property, it is not required that the
 common hazard rate of $Y$'s is  necessarily $\bar \lmd$ in the above theorems. In fact,  the $Y$'s could be a random sample from any continuous distribution. }
 \end{remark}

 The above theorem can be extended to the PHR model using the technique used in {\em Genest, Kochar and Xu} (2009).

 \begin{theorem}\label{dep3}\rm \rm Let $X_1,\ldots,X_n$ be independent continuous random variables following the
PHR model with $\lmd_1, \ldots, \lmd_n$ as the proportionality parameters. Let $Y_1,\ldots,Y_n$ be i.i.d. continuous random
variables with common survival function $\olf^{\bar\lmd}$, then
\be
(X_{i:n}|X_{1:n})\prec_{\rm SI}(Y_{i:n}|Y_{1:n}) \, \mbox{ for $2\le i \le n$}.\ee
\end{theorem}

\NI{\bf Proof.} \, Let $R(t)=  - \log \{ {\bar F} (t) \}$ in the PHR model.  Make the monotone transformations $X_i^* = R(X_i)$ and  $Y_i^* =  - \log \{ {\bar F} (Y_i) \}$.
Then the transformed  variable $X^*_i$ has exponential distribution  with hazard rate $\lmd_i$, \, $i=1,\ldots,n$ and $Y_1^*, \ldots, Y_i^*$ is a random sample from an exponential distribution with parameter ${\bar \lmd}$. Let  $X^*_{(1)} < \cdots < X^*_{(n)}$ and $Y^*_{(1)} < \cdots <
Y^*_{(n)}$ be the order statistics corresponding to the new sets of variables.

In view of their invariance by monotone increasing transformations of the margins, the copulas associated with the pairs $(X_{1:\,n}, X_{n:\,n})$ and $(X^*_{1:\,n}, X^*_{n:\,n})$ coincide. Similarly, the pairs $(Y_{1:\, n}, Y_{n:\, n})$ and $(Y^*_{1:\,n}, Y^*_{i:\,n})$ have the same copula.

Also, since  the more SI dependence ordering is copula-based,
$$
(X_{n:n}|X_{1:n})\prec_{\rm SI}(Y_{n:n}|Y_{1:n})  \Leftrightarrow  (X^*_{n:n}|X^*_{1:n})\prec_{\rm SI}(Y^*_{n:n}|Y^*_{1:n})
$$
\qed

Under the conditions of Theorem \ref{dep3}, an upper bound on $\kappa (X_{1:n}, X_{(i:n})$ is given by $\kappa (Y_{1:n}, Y_{i:n})$.  {\em Av\'{e}rous, Genest and Kochar} (2005) obtained an analytic expression for computing the exact values of the  Kendall's  $\tau$ for any pair of order statistics from a random sample from a continuous distribution which in our case reduces to
\be
\tau(Y_{1:n}, Y_{i:n}) = 1-\frac{ 2(n-1)}{2n-1}\left({n-2}\atop {i-2}\right)
\sum_{s=0}^{n-i}\left({n} \atop {s}\right)  {\Big /}
\left( {2n-2} \atop {n-i+s} \right).
\ee

\end{document}